\documentclass[12pt]{article}
\usepackage[final]{epsfig}	
\usepackage{amsfonts,amssymb,latexsym,
amscd, theorem,epsfig,graphics}
\newtheorem{theorem}{Theorem}
\newtheorem{lemma}{Lemma}[section]

\newtheorem{remark}[lemma]{Remark}

\begin{document}
\newcommand{\eps}{{\varepsilon}}
\newcommand{\proofend}{$\Box$\bigskip}
\newcommand{\C}{{\mathbf C}}
\newcommand{\Q}{{\mathbf Q}}
\newcommand{\R}{{\mathbf R}}
\newcommand{\Z}{{\mathbf Z}}
\newcommand{\RP}{{\mathbf {RP}}}

\title {On skew loops, skew branes and quadratic hypersurfaces}
\author{Serge Tabachnikov \thanks{ Supported, in part, by a BSF grant}\\
{\it Department of Mathematics, Penn State University}\\
{\it University Park, PA 16802, USA}\\
e-mail: {\it tabachni@math.psu.edu} 
}
\date{}
\maketitle

\begin{abstract} A skew brane is an immersed codimension 2 submanifold in affine space,
free from pairs of parallel tangent spaces. Using Morse theory, we prove that a skew brane
cannot lie on a quadratic hypersurface. We also prove that there are no skew loops on
embedded ruled developable discs in 3-space.

{\it MSC}: 53A05, 53C50, 58E05

{\it Key words}: skew loops and skew branes, quadratic hypersurfaces, double normals, Morse
theory, developable surfaces
\end{abstract}
\bigskip

\section{Introduction and formulation of results}

In a lecture, given in 1966,  H. Steinhaus asked the question: Do there exist smooth closed
space curves without a pair of parallel tangents? Call such a curve a {\it skew loop}. 
Steinhaus conjectured that skew loops did not exist. The problem was solved by B. Segre in
\cite{Se 1, Se 2}. Segre constructed examples of skew loops in 3-space but proved that no
skew loops lie on an ellipsoid or an elliptic paraboloid. Thus skew loops have an  aversion
to quadrics.

These results were recently extended and refined by M. Ghomi and B. Solomon  \cite{G-S}. On
the one hand, the following theorem holds: {\it Nonnegatively curved quadrics admit no
 skew loops.} In particular, no skew loops lie on  one sheet of a 2-sheeted hyperboloid.
On the other hand, the converse is true: {\it Any non-quadratic surface with a point of
positive curvature admits a skew loop.} 

Another recent result on this subject  \cite{Wu}: {\it Every knot type can be realized by a
skew loop}.

A multi-dimensional version of the problem reads as follows. Let $f: M^{n} \to V^{n+2}$ be
an immersion of a manifold to an affine space. One asks whether there are pairs of points
$x,y \in M$ such that the tangent spaces $df (T_{x} M)$ and
$df (T_{y} M)$ are parallel in $V$. Such pairs correspond to self-intersections of the
image of the tangent Gauss map $M \to G_{n} (V)$ in the Grassman manifold $G_n (V)$ of
$n$-dimensional non-oriented subspaces in $V$. Since dim  $G_n (V) = 2n$, these
self-intersections generically occur in isolated points. An immersed closed connected
codimension 2 manifold is called a {\it skew brane} if it is free from pairs of parallel
tangents.  

J. White \cite {Wh} considered the case of $M^n$ immersed into the unit sphere
$S^{n+1}$ in Euclidean space $V^{n+2}$. The main observation in this case is as follows: if
$xy$ is a double normal of $M$, that is, the segment $xy$ is perpendicular to $M$ at both
end points, and the points $x$ and $y$ are not antipodal then the tangent spaces
$T_x M$ and $T_y M$ are parallel in $V$. This made it possible to estimate below the
number of parallel tangents of $M$ by Morse theory of the distance squared function on $M
\times M$, studied in a related paper \cite{T-W}. See \cite{Pu 1, Pu 2} for  recent results
on double normals of immersed submanifolds in Euclidean space.

The goal of this note is to extend the above mentioned results on non-existence of skew
loops and skew branes from convex quadratic surfaces to  quadratic hypersurfaces of all
signatures. Here is our main result.

\begin{theorem} \label{main}  Let $Q^{n+1}$ be a non-degenerate quadratic hypersurface in a
vector space $V^{n+2}$ and $f: M^{n} \to Q$ be a generic smooth immersion of a closed
manifold. Then there exist distinct points $x,y \in M$ such that the tangent spaces $df
(T_{x} M)$ and $df (T_{y} M)$ are parallel in $V$.
\end{theorem}

In particular, there are no generic skew loops on  one-sheeted hyperboloids and hyperbolic
paraboloids, a new result compared  to  \cite {G-S, Se 1, Se 2}. The meaning of the general
position assumption is explained in section 2.

We give applications to two types of quadratic surfaces in 3-space: to one sheet of the 
2-sheeted hyperboloid $z^2-x^2-y^2=1, z>0$ and to 1-sheeted hyperboloid $x^2+y^2-z^2=1$.

\begin{theorem} \label{hyper} {\rm (i)} Let $\gamma$ be a generic smooth closed immersed
curve on one sheet of the 2-sheeted hyperboloid, and let $a$ be the number of its double
points.  Then $\gamma$ has at least $a+2$ pairs of  tangent lines, parallel in the ambient
space. \break  {\rm (ii)}  Let $\gamma$ be a generic smooth closed immersed curve on the
1-sheeted hyperboloid, and assume that $\gamma$ is space-like: the quadratic form
$x^2+y^2-z^2$, evaluated on its tangent vector, is everywhere positive. Let $a$ be the
number of double points of
$\gamma$ and $b$ the number of intersections with its antipodal image. Then  $\gamma$ has 
at least $a+b$ pairs of  tangent lines, parallel in the ambient space. If $a=b=0$ then
there are at least 2 pairs of parallel tangents. 
\end{theorem}

Statement (i) of Theorem \ref{hyper} holds for spherical curves that do not intersect their
antipodal images -- see \cite{Wh}.

Another result concerns non-existence of skew loops on developable surfaces in
3-dimensional Euclidean space. Let $D^2$ be a disc, embedded in the Euclidean space
$V^3$ as a smooth ruled developable surface, and let $\gamma$ be a generic smooth closed
immersed curve in $D$. 

\begin{theorem} \label{devel} There exists a pair of points on $\gamma$ with tangent lines,
parallel in the ambient space.
\end{theorem}

Since cylinders are ruled developable surfaces, this theorem is in sharp contrast with a
result from \cite{G-S}: {\it the cylinder over any non centrally symmetric, strictly convex
plane oval contains a skew loop}. This skew loop goes around the cylinder and is not
contained in a simply connected domain. Moreover, it was pointed out by D. Burago that Theorem
\ref{devel} fails if the the developable surface contains flat pieces; in this case it
may not be ruled.
\smallskip

{\bf Acknowledgments}. It is a pleasure to acknowledge the interest and help of D. Burago,
S. Duzhin, M. Ghomi, J. Landsberg, A. Malyutin, R. Matveyev, P. Pushkar' and especially B.
Solomon who introduced me to the subject. Part of work was done while I was visiting
Max-Planck-Institut fur Mathematik in Bonn; I am grateful to the Institute for its
invariable hospitality.

\section{Proof of Theorem \ref{main}}

Assume that 
$$Q^{n+1} = \{x \in V^{n+2}|\ x \cdot x = 1\}$$   
where dot denotes a (pseudo)Euclidean structure of signature $(p+1,q),\ p+q = n+1,\ p,q \geq
0$. The quadratic hypersurface $Q$ is the unit (pseudo)sphere. 
 Let $(x_1,\xi_1)$ and $(x_2,\xi_2)$ be two contact elements of $Q$; here 
$x_i \in Q$ is the foot point and $\xi_i \subset T_{x_i} Q$ is a hyperplane in the tangent
space to $Q$.  Assume that $x_2 \neq \pm x_1$, and let $U$ be the 2-dimensional subspace
spanned by $x_1$ and $x_2$.

\begin{lemma} \label{par} The spaces $\xi_1$ and $\xi_2$ are parallel in $V$ if and only if
they are both perpendicular to $U$. 
\end{lemma}

{\bf Proof}. Two $n$-dimensional spaces with the same orthogonal complement  are parallel.
Conversely, let $\xi_1$ and $\xi_2$ be parallel. Then their orthogonal complements in $V$
coincide. Note that the position vector $x_i$ is  orthogonal to the (pseudo)sphere $Q$ and,
therefore, to $\xi_i$. Hence the orthogonal complement of
$\xi_i$ contains $x_i$, and the common orthogonal complement of $\xi_1$ and $\xi_2$ is
$U$.
\proofend

\begin{remark} \label{geo} {\rm One can reformulate Lemma \ref{par} as follows. Consider
the 2-parameter family of parallel codimension 2 affine subspaces in $V$, orthogonal to a
given 2-dimensional vector space $U$. Then the locus of tangency of these affine subspaces
with $Q$ is the intersection curve $U \cap Q$. This curve is a geodesic on
$Q$ for the (pseudo)metric, induced from $V$.}
\end{remark}

The next lemma extends a result of White \cite{Wh} from spheres to all non-degenerate
quadratic hypersurfaces. Let $f: M^{n} \to Q$ be an immersion. Consider the function on $M
\times M$ given by the formula: 
\begin{eqnarray}
\phi(x,y) = f(x) \cdot f(y). \label{phi} 
\end{eqnarray}

\begin{lemma} \label{crit} Let $f(y)  \neq \pm f(x)$. The spaces $df (T_{x} M)$ and $df
(T_{y} M)$ are parallel in $V$ if and only if  $(x,y)$ is a critical point of the function
$\phi$.
\end{lemma}

{\bf Proof}. Let $U$ be the the 2-dimensional space spanned by $f(x)$ and $f(y)$. A point 
$(x,y)$ is  critical for the function $\phi$ if and only if $f(x)$ is orthogonal to $df
(T_{y} M)$ and $f(y)$ is orthogonal to $df (T_{x} M)$. Also
$f(x)$ is orthogonal to $T_{f(x)} Q$ and hence to $df (T_{x} M)$, and likewise, $f(y)$ is
orthogonal to $df (T_{y} M)$. Thus $(x,y)$ is  critical for $\phi$ if and only if $U$ is
orthogonal to the spaces $df (T_{x} M)$ and $df (T_{y} M)$.  By Lemma \ref{par}, this holds
if and only if these two spaces are parallel.
\proofend

\begin{remark} \label{graph} {\rm A special degenerate case of a quadratic hypersurface is the
graph  of a quadratic form. Let $V = W \oplus \R$ where space $W$ has a non-degenerate
quadratic form $\tau$, and let 
$$Q = \{(x,z) \in V|\ z = \tau(x)\}.$$
 Denote by $\pi$ the projection of $V$ to $W$. Let $f: M^{n} \to Q$ be an immersion and let 
$g=\pi \circ f$. Consider the function on $M \times M$ given by the formula: 
$$\psi(x,y) = \tau(g(x)-g(y)).$$
 Then an analog of Lemma \ref{crit} reads as follows: for $g(y)  \neq  g(x)$,  the spaces
$df (T_{x} M)$ and $df (T_{y} M)$ are parallel in $V$ if and only if  $(x,y)$ is a critical
point of the function $\psi$.}
\end{remark}
 
The proof of Theorem \ref{main} proceeds as follows. Let $f: M^{n} \to Q^{n+1}$ be an
immersed skew brane such that the self-intersection of $M$ and its intersection with the
antipodal image are transversal. This is the desired general position property. Consider the
double and the antipodal double loci in
$M
\times M$:
$$A = \{ (x,y)\in M \times M|\ y \neq x,\ f(y) = f(x) \},$$ and
$$B = \{ (x,y)\in M \times M|\ f(y) = -f(x) \},$$ and let $D \subset M \times M$ be the
diagonal. Then $A$ and $B$ are closed submanifolds of dimensions $n-1$ (cf. \cite{Wh}).
Furthermore we require that $\phi$ in (\ref{phi}) is a Morse function on $M \times M - A
\cup B \cup D$; this is also provided  by the genericity of the immersion.

By Lemma \ref{crit}, the function $\phi$ has no critical points on $M
\times M$ off the sets $A, B$ and $D$. We will show that this cannot be the case, and this
contradiction will prove Theorem \ref{main}.

Recall that if a smooth function has a non-degenerate critical manifold then the tangent
space to the critical manifold is the kernel of its Hessian. By  Morse index  one means the
index of the Hessian in the normal bundle to the critical manifold -- see, e.g., \cite{Bo}. 

\begin{lemma} \label{ind} The manifolds $D, A$ and $B$ are critical manifolds of the
function $\phi$ with the critical values $1, 1$ and $-1$, respectively. The manifolds
$A$ and $B$ are non-degenerate, and the respective Morse indices equal $p$ and $q$.  Let
$(x,x)$ be a point of the diagonal $D$, and assume that the  (pseudo)Euclidean structure on
$T_x M$, induced by $f$ from $V$,  is non-degenerate and has index
$r$ (which can be equal to either $q$ or $q-1$). Then  the Morse index of $D$ at point
$(x,x)$ equals $n-r$.

\end{lemma}

{\bf Proof}. Consider a point $(x_1, x_2) \in M \times M$, and let $(u_1, u_2)$ be a
tangent vector to $M \times M$ at this point. Let $(\gamma_1 (t),\gamma_2 (t))$ be a curve
on $M \times M$ such that $\gamma_i(0)=x_i$ and $\gamma'_i(0)=u_i; \ i=1,2$. Set:
$q_i(t) = f(\gamma_i(t))$, and introduce the following notation: $q_i(0)=y_i, q'_i(0)=v_i$
and $q''_i(0)=\eta_i$. In particular, $f(x_i)=y_i$ and $df(u_i)=v_i$.

To compute the Hessian, we make calculations modulo  $t^3$. One has:
$$q_i(t) = y_i + t v_i + {{1}\over{2}} t^2 \eta_i,$$ and therefore 
$$
\phi(\gamma_1(t), \gamma_2(t)) =  y_1 \cdot y_2 + t (v_1 \cdot y_2 +  y_1 \cdot v_2)
+{{t^2}\over{2}} (\eta_1 \cdot y_2 + 2 v_1 \cdot v_2 +  y_1 \cdot \eta_2).
$$  If $(x_1, x_2) \in D \cup A \cup B$ then $y_2 = \pm y_1$ (plus for $D \cup A$ and minus
for $B$), and the linear terms vanish:
$$v_1 \cdot y_2=y_1 \cdot v_2=0.$$ Hence $D, A$ and $B$ are critical manifolds. Since
$q_i(t) \cdot q_i(t) = 1$ one has:
$$y_i  \cdot \eta_i + v_i \cdot v_i =0.$$ Therefore 
\begin{eqnarray}
\phi(\gamma_1(t), \gamma_2(t)) = 1-{{t^2}\over{2}} (v_1 -  v_2) \cdot (v_1 -  v_2)
\label{ind1}
\end{eqnarray} for  $(x_1, x_2) \in D \cup A$ and 
\begin{eqnarray}
\phi(\gamma_1(t), \gamma_2(t)) = -1 + {{t^2}\over{2}} (v_1 + v_2) \cdot (v_1 + v_2) 
\label{ind2}
\end{eqnarray} for  $(x_1, x_2) \in B$.

Assume that $(x_1, x_2) \in A$. According to (\ref{ind1}), the Hessian of $\phi$ is the
quadratic form on $T_{x_1} M \oplus T_{x_2} M$ which is the composition
$$T_{x_1} M \oplus T_{x_2} M \to T_{y_1} Q \to \R,$$ where the first map is a linear
epimorphism  given by the formula 
\begin{eqnarray} (u_1, u_2) \to df (u_1) - df(u_2), \label{map}
\end{eqnarray}
 and the second map is negative the quadratic form on $Q$ given by the restriction of the
(pseudo)Euclidean structure in $V$. The kernel of (\ref{map}) is the tangent space to $A$
at point $(x_1, x_2)$, and the Morse index is equal to the index of  the restriction of the
(pseudo)Euclidean structure to $Q$. This structure in $V$ has  signature $(p+1,q)$, and $Q$
is given by the equation $x \cdot x =1$. Therefore the restriction of the (pseudo)Euclidean
structure to $Q$ has signature
$(p,q)$, and the Morse index is equal to $p$. 

A similar argument, with (\ref{ind2}) replacing (\ref{ind1}), shows that the Morse index of
the critical manifold $B$ equals $q$. 

An analogous argument works for the diagonal as well:  by (\ref{ind1}), the Hessian of
$\phi$ is the composition
$$T_{x} M \oplus T_{x} M \to T_{f(x)} Q \to \R,$$ where the first arrow (\ref{map}) is a
linear map onto $df (T_{x} M) \subset T_{f(x)} Q$. The restriction of the quadratic form
from $Q$ to $df (T_{x} M)$ has signature either
$(p-1,q)$ or $(p,q-1)$, depending on the direction (``space-like" or ``time-like") of $M$
at point $f(x)$. Thus $r$ is either $q$ or $q-1$.
\proofend

We are ready to complete the proof of Theorem \ref{main}. Recall that we work toward a
contradiction under the assumption that the function $\phi$ has no critical points outside
the manifolds $A, B$ and $D$. Consider two cases. 

Case 1: $q>0$. Our assumption implies that $\phi$ attains minimum $-1$ on $B$. On the other
hand, by Lemma \ref{ind}, the Morse index of $B$ equals $q$, a positive number, hence $B$
is not a minimal critical manifold.

Case 2: $q=0, p=n+1$. Then the quadratic hypersurface $Q$ is the unit sphere, the cases
studied in \cite{Wh}. Note that, by Lemma \ref{ind}, the diagonal $D$ is a non-degenerate
critical manifold with the Morse index $n$. 

Recall the Morse inequalities for Morse-Bott functions, that is, functions with
non-degenerate critical manifolds \cite{Bo}. Let $g$ be such a function on a closed
manifold $L$ with critical submanifolds $N_i,\ i=1, \dots, k$, whose Morse indices are equal
to $\mu_i$.  Let $P_t$ denote the Poincar\'e polynomial in variable
$t$ of a manifold with coefficients in a fixed field. Then the Morse inequalities can be
written as follows:
\begin{eqnarray}
 \sum_{i=1}^k t^{\mu_i} P_t (N_i) = P_t (L) + (1+t) Q_t \label{ineq}
\end{eqnarray} where $Q_t$ is a polynomial with non-negative coefficients.

Apply (\ref{ineq}) in our situation. Let
$$P_t(A) = a_0 + a_1 t + \dots + a_{n-1} t^{n-1}, P_t(B) = b_0 + b_1 t + \dots + b_{n-1}
t^{n-1}$$ and
$$P_t(D) = P_t(M) = 1 + d_1 t + \dots + t^n.$$ Then the Morse inequalities read:
$$ t^n (1 + d_1 t + \dots + t^n) + t^{n+1} (a_0 + a_1 t + \dots + a_{n-1}) + (b_0 + b_1 t +
\dots + b_{n-1} t^{n-1}) =$$
$$ (1 + d_1 t + \dots + t^n)^2 + (1+t) Q_t. $$ The coefficient of $t^n$ on the left hand
side is 1, while on the right hand side it is not less than 2. This is a desired
contradiction proving Theorem \ref{main}.

\section{Loops on surfaces. Open questions}

Let us prove Theorem \ref{hyper}. In the case of the hyperboloid $z^2-x^2-y^2=1$, one has:
$p=0, q=2$. The set $A$ consists of $2a$ isolated points (a self-intersection counts twice:
as $(x,y)$ and as $(y,x)$), each with the Morse index 0, and the set $B$ is empty. The
Morse index of the diagonal $D$ equals 0. Let $2d_0, 2d_1$ and $2d_2$ be the number of
critical points of the function $\phi$ outside $D \cup A$; these numbers are even since
$\phi$ is invariant under the involution $(x,y) \to (y,x)$. The Morse inequalities read:
$$(1+t) + 2a + (2d_0 + 2d_1t + 2d_2 t^2) = (1+t)^2 + (1+t)(q_0 + q_1 t)$$ where $q_0, q_1
\geq 0$. It follows that 
$$q_0 = 2a + 2d_0 \geq 2a,\ 2d_1 = 1+ q_0 +q_1\geq 1+2a,\ 2d_2
\geq 1.$$
 Hence $d_1 \geq 1+a,\ d_2 \geq 1$, and $d_1 + d_2 \geq 2+a$, as claimed.

Consider the hyperboloid $x^2+y^2-z^2=1$. In this case $p=q=1$. The sets $A$ and
$B$ consist of $2a$ and $2b$ isolated points, respectively, each with the Morse index 1. If
$\gamma$ is space-like then the Morse index of the diagonal $D$ equals 0. The Morse
inequalities read:
$$(1+t) + (2a+2b)t + (2d_0 + 2d_1t + 2d_2 t^2) = (1+t)^2 + (1+t)(q_0 + q_1 t)$$ where $q_0,
q_1 \geq 0$. It follows that
$$2d_0 = q_0,\ 2a+2b+2d_1 = 1+q_0+q_1,\ 2d_2=1+q_1.$$ If $a=b=0$ then $d_1 \geq 1,\ d_2
\geq 1$, and $d_1 + d_2 \geq 2$. Also $2d_0 + 2d_2 = 1+q_0+q_1 \geq 2a+2b$, and thus $d_0 +
d_2 \geq a+b$.
\smallskip

Now we proceed to the proof of Theorem \ref{devel}. Recall the classification of ruled 
developable surfaces (see, e.g., \cite{St}). Such a surface is either a cylinder or a cone
or consists of tangent lines to a space curve (called the edge of regression). In either
case,  the tangent space to the surface remains the same along each ruling.

Let $D$ be a ruled embedded developable disc and $\gamma$ a closed immersed curve in it.  Then
$D$ is foliated by straight lines, the ruling of the developable surface. Unfold $D$ to the
plane to obtain a plane domain $\bar D$ foliated by straight lines and an immersed curve
$\bar \gamma$ in it. Assume that, for two distinct intersection points $x, y$ of
$\bar \gamma$ with a leaf $\bar l$ of the foliation, the tangents $T_x \bar \gamma$ and
$T_y \bar \gamma$ are parallel. Then the respective tangent lines to $\gamma$ lie in the
plane, tangent to $D$ along the corresponding ruling $l$, and are parallel in this plane.
Hence they are parallel in the ambient space. Thus  Theorem \ref{devel} follows from the
next lemma.

\begin{lemma} \label{fol} Let $\bar \gamma$ be a generic closed immersed curve in a  simply
connected plane domain, foliated by straight lines. Then there exist two distinct points
$x, y \in \bar \gamma$ on the same leaf of the foliation, such that the tangents
$T_x \bar \gamma$ and $T_y \bar \gamma$ are parallel.
\end{lemma}

{\bf Proof}. This proof was proposed by A. Malyutin; another proof was given by R.
Matveyev.

\begin{figure}[ht]
\centerline{\epsfbox{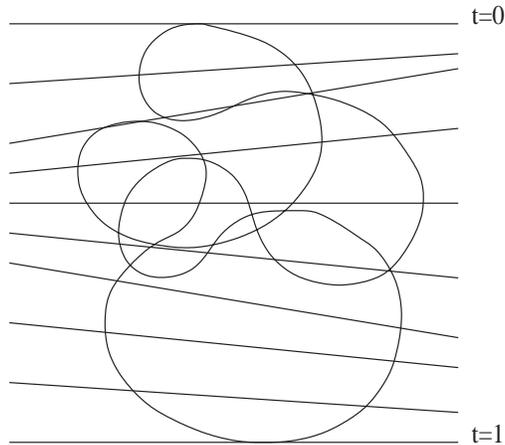}}
\caption{A curve in a domain foliated by lines}
\end{figure}

\begin{figure}[ht]
\centerline{\epsfbox{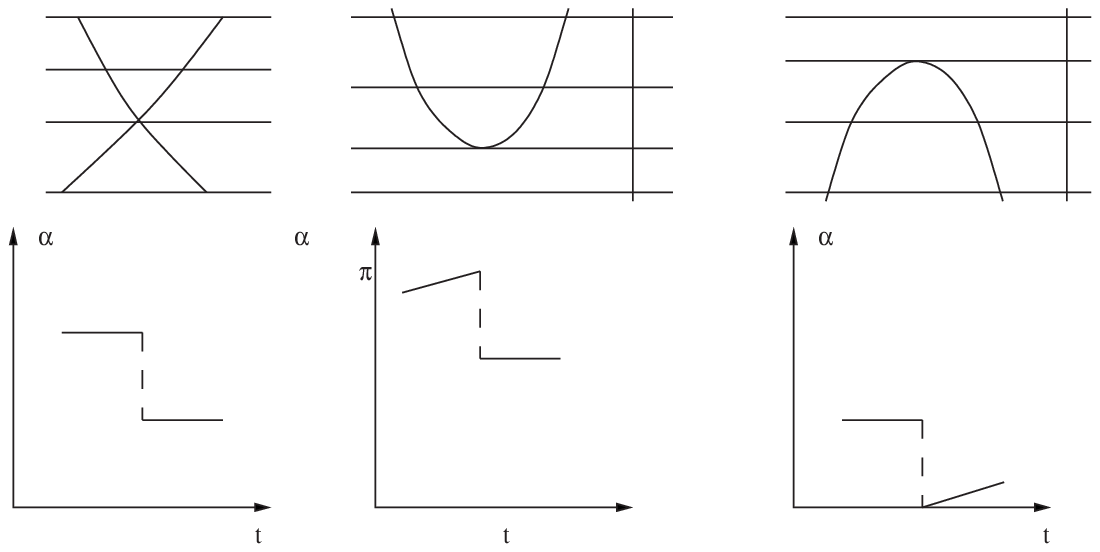}}
\caption{Discontinuities of the functions $\alpha_{\pm} (t)$}
\end{figure}

Let $\bar l_t,\ t \in [0,1]$, be the family of lines foliating the domain;  we think of
these lines as oriented from left to right. Assume that $\bar \gamma$ lies between
$\bar l_0$ and $\bar l_1$, and that these lines touch the curve -- see Figure 1. Let $C_{+}
(t)$ and $C_{-} (t)$ be the rightmost and the leftmost points of $\bar \gamma \cap \bar
l_t$. The curves $C_{\pm} (t)$ are piecewise smooth. Let $\alpha_{\pm} (t) \in [0,\pi]$ be
the angle between $C_{\pm} (t)$ and $\bar l_t$.

The functions $\alpha_{\pm} (t)$ are not continuous but their discontinuities are easily
described. The three types of discontinuities of $\alpha_-(t)$ are shown in Figure 2: in
each case the graph of  $\alpha_-(t)$ has a descending vertical segment. Likewise, there
are three types of discontinuities of $\alpha_+(t)$, and in each case the graph has an
ascending vertical segment.

Note that, for $t$ close to 0,  $\alpha_-(t)$ is  close to 0, while 
$\alpha_+(t)$ is  close to $\pi$. Similarly, for $t$ close to 1, 
$\alpha_-(t)$ is  close to $\pi$ and  $\alpha_+(t)$ is  close to 0. 

We claim that that the graphs of $\alpha_-(t)$ and $\alpha_+(t)$ have a common point not on
a vertical segment of either graph. Approximate both functions by smooth ones so that the
vertical segments of the graphs become very steep: $\alpha'_-(t) \ll 0$ and 
$\alpha'_+(t) \gg 0$ on the corresponding small intervals. Consider the  function
$\beta(t) = \alpha_+(t) - \alpha_-(t)$. For $t$ near 0, $\beta(t) >0$, and for
$t$ near 1, $\beta(t) < 0$. Let $t_0$ be the smallest zero of $\beta(t)$. Since $\beta$
changes sign from positive to negative, $\beta'(t_0)
\leq 0$. Therefore $t_0$ is not on the segments where $\alpha'_-(t) \ll 0$ or $\alpha'_+(t)
\gg 0$. Thus $t_0$ is a desired common value of the two angles.
\proofend

\begin{remark} \label{var} {\rm If the foliation consists of parallel lines, that is,  the
developable surface is (part of) a cylinder one can prove Lemma \ref{fol}  as follows. Let
$xy$ be the longest segment of a leaf bounded by its intersection with the curve $\bar
\gamma$. Then the tangents $T_x \bar \gamma$ and $T_y
\bar \gamma$ are parallel.

Likewise, assume that the foliation consists of the lines passing through a point $0$, that
is,  the developable surface is (part of) a cone. Let
$xy$ be the segment of a leaf, $x,y \in \bar \gamma$, such that the ratio $|Ox|/|Oy|$ is
greatest possible. Then the tangents $T_x \bar \gamma$ and $T_y \bar \gamma$ are parallel.}
\end{remark}

\begin{remark} \label{Bur} {\rm Let us describe an example of a simple skew loop on a
developable disc, suggested by D. Burago. This developable surface contains flat regions and
does not admit a continuous foliation by straight segments. 

Take an equilateral triangle in the plane; its interior will be
folded to make the desired embedding. Draw three lines near each corner, cutting
across and not parallel to opposite sides. Fold smoothly along these  lines so that the
obtained surface consists of a flat hexagon with 3 triangular flaps (fins?) going
more-or-less vertically. The curve $\gamma$ is a simple smooth convex curve that goes
along the perimeter of the original triangle. After the flaps are folded, the curve doesn't
have parallel tangents in the ambient space.}
\end{remark}

We now pose some problems on skew loops and skew branes.
\smallskip

{\bf Problem 1}. Give lower bounds on the number of pairs of parallel tangent spaces of $M$
in Theorem \ref{main} in terms of the topology of $M$. Such lower bounds on the number of
double normals of a generic immersed submanifold in Euclidean space are found in \cite
{T-W, Pu 1, Pu 2} (the estimates by Pushkar' are sharp in some cases).
\smallskip

{\bf Problem 2}. Extend Theorem \ref{main}  and lower bounds of Problem 2 to the case when
$M$ is a wave front, Legendrian isotopic to a smooth manifold. We refer to \cite {Pu 1, Pu
2} for such results on double normals. 
\smallskip

{\bf Problem 3}. Our proof of Theorem \ref{main} is based on Morse theory of the distance
function $\phi: M \times M \to \R$. The proofs in \cite {Pu 1, Pu 2} make use of the
function $\psi: M \times M \times S^N \to \R$ given by $\psi(x,y,v) = (x-y) \cdot v$; here
$M$ is immersed in Euclidean space $\R^{N+1}$. Critical points of $\psi$ are double
normals, and $\psi$ is invariant under the free action of $\Z_2$ on the manifold $M \times
M \times S^N$ sending $(x,y,v)$ to $(y,x,-v)$.  Can one apply this approach in the set-up
of Theorem \ref{main}?
\smallskip

{\bf Problem 4}. An immersion $f: M \to Q$ determines an immersion of $M$ to the space of
geodesics $G$ of the pseudo-metric on $Q$: one assigns to a point of $M$ its normal line.
One is looking for self-intersections of the image of $M$ in
$G$. If the pseudo-metric on $Q$ is a metric (then $Q$ is either the sphere or the
hyperbolic space) then the space of oriented geodesics $\tilde G$ has a symplectic
structure, and  $M$ is immersed to $\tilde G$ as a Lagrangian submanifold. What can be said
in the pseudo-metric case?
\smallskip

{\bf Problem 5}. Construct examples of skew branes. For instance, let $N^n \subset
W^{n+1}$  be a strictly convex {\it non centrally symmetric} closed hypersurface in vector
space, and let $n$ be odd. Can one construct a smooth function $g: N \to \R$ such that its
graph in $V = W \times \R$ is a skew brane? For $n=1$ this is the Cylinder Lemma of
\cite{G-S}. 

{\bf Problem 6}. Do Theorems \ref{main} and \ref{devel} hold without the general position
assumption (conjecturally, they do)?  The property of being a skew loop or a skew brane is
not open: a small perturbation of a skew loop may create  parallel tangents at a pair of close
points $(x,y) \in M \times M$ near the diagonal $D \subset M \times M$. Thus one cannot simply
argue by deforming a skew loop or a skew brane to a generic one.

\end{document}